\documentclass[12pt, a4paper]{article}
\usepackage[latin1]{inputenc}
\usepackage[english]{babel}
\usepackage{amsmath}
\usepackage[english]{babel}
\usepackage{indentfirst}
\usepackage{amstext}
\usepackage{enumerate}
\usepackage{amsfonts}
\usepackage{textcomp}
\usepackage{amssymb}
\usepackage[dvips]{graphicx}
\usepackage{setspace}
\usepackage[all]{xy}

\newcommand {\demo}{\hskip -0.6cm{\bf Proof.  }}
\newcommand {\fim}{\hfill{$\square$}\vskip 1pc}
\newcommand {\nl}{\newline}

\newcommand {\N}{\mathbb{N}}
\newcommand {\Z}{\mathbb{Z}}

\newcommand {\F}{\mathbb{F}}

\newcommand {\GG}{\mathcal{G}}

\newtheorem{teorema}{Theorem}[section]
\newtheorem{lema}[teorema]{Lemma}
\newtheorem{corolario}[teorema]{Corollary}
\newtheorem{definicao}[teorema]{Definition}
\newtheorem{proposicao}[teorema]{Proposition}
\newtheorem{exemplo}[teorema]{Example}
\newtheorem{remark}[teorema]{Remark}
\begin{document}
\onehalfspace

\title{Infinite alphabet edge shift spaces via ultragraphs and their C*-algebras}
\maketitle
\begin{center}

{\large Daniel Gon\c{c}alves\footnote{This author is partially supported by CNPq.} and Danilo Royer}\\
\end{center}

\vspace{8mm}

\abstract We define a notion of (one-sided) edge shift spaces associated to ultragraphs. In the finite case our notion coincides with the edge shift space of a graph. In general, we show that our space is metrizable and has a countable basis of clopen sets. We show that for a large class of ultragraphs the basis elements of the topology are compact. We examine shift morphisms between these shift spaces, and, for the locally compact case, show that if two (possibly infinite) ultragraphs have edge shifts that are conjugate, via a conjugacy that preserves length, then the associated ultragraph C*-algebras are isomorphic. To prove this last result we realize the relevant ultragraph C*-algebras as partial crossed products.

\vspace{1.0pc}

\vspace{1.0pc}

\section{Introduction}

In classical symbolic dynamics a shift space is a set of infinite words that represent the evolution of a discrete system: One starts with a finite alphabet $A$, with the discrete topology, and constructs the infinite products $A^\Z$ and $A^\N$. The map $\sigma$ attached to these spaces shifts all the entries of the sequences one to the left. A {\em shift space} (or a {\em subshift}) is then a closed subspace of $A^\Z$ or $A^\N$ which is invariant under $\sigma$. These dynamical systems are well-studied -- see \cite{LindMarcus} for an excellent reference.

A possible approach that can be used to define analogues of shift spaces to infinite countable alphabets is to consider $A^\N$ with the product topology, which yields a space that is not locally compact. This difference in behavior, when compared to the finite alphabet case, has motivated numerous efforts, by various authors, in the search of a definition (and consequent study) of an analogue of shift spaces when dealing with infinite countable alphabets. Examples of these approaches are \cite{FiebigFiebig1995}, where the authors use the Alexandroff one-point compactification for locally compact shift spaces over countable alphabets (such compactification was used to prove several results for the entropy of countable Markov shifts \cite{Fiebig2001,Fiebig2003,FiebigFiebig1995,FiebigFiebig2005}) and \cite{OTW}, where the authors use the Alexandroff compactification of the alphabet to obtain a new shift space (the theory of these shift spaces were further developed in \cite{GRUltra, Mstep, GSS2, GSS, GSS1}).

Turning back to the finite alphabet case for a minute, recall that a shift space can be defined in terms of a set of forbidden words. 
If the set of forbidden words is finite the shift is called a shift of finite type. Shifts of finite type are among the most important in symbolic dynamics, having practical applications such as finding efficient coding schemes to store data on computer disks, see \cite{LindMarcus}. It is well known that a shift space is a shift of finite type if, and only if, it is conjugate to the edge shift coming from a finite graph with no sinks if, and only if, it is conjugate to a 1-step shift (which can be seen as shifts associated to a matrix of 0-1 and are also called Markov shifts). Taking the above in consideration, it is not surprising that analogues of shift spaces to infinite alphabets appeared most commonly in the context of countable-state Markov chains (or equivalently, shifts coming from countable directed graphs or matrices). Examples of such an approach can be found in \cite{FF02, FiebigFiebig2005, kitchens1997}. Further development of the theory appeared in \cite{BBG06, BBG07} where Boyle, Buzzi and G\'omez study almost isomorphism for countable-state Markov shifts, \cite{CS09, IY12, MU01, Sa99} where thermodynamical formalism for such shifts is developed, \cite{GS} where sliding block codes are characterized, and many other papers (see \cite{OTW} for a comprehensive list). 

It is our point of view that an analogue of a shift of finite type should be connected with C*-algebras, in a way similar to how Markov shifts are connected to Cuntz-Krieger algebras in the finite alphabet case. In this line of thought Exel and Laca, see \cite{ExelLaca}, define a C*-algebra $\mathcal{O}_B$ from a given countable $\{0,1\}$-matrix $B$ which is thought of as the incidence matrix of an infinite graph, and propose that the spectrum of a certain commutative C*-subalgebra of $\mathcal{O}_B$ is a good candidate for the Markov shift associated to the graph. In \cite{OTW}, motivated by work of Paterson and Welch (\cite{PW}), Ott-Tomforde-Willis introduce one sided shift spaces associated to infinite alphabets and connect then to graph C*-algebras. In the case of a directed graph their space is related, but is not the same, as the path space studied by Webster, see \cite{Webster}. In fact, the boundary path space studied by Webster is the spectrum of a certain commutative C*-subalgebra of the graph C*-algebra. 

The above considerations led us to believe that, in connection with C*-algebra theory, the analogue of a shift of finite type should be the spectrum of a commutative subalgebra of a combinatorial C*-algebra. In this context the notion of ultragraph C*-algebras, which were defined by Tomforde (see \cite{Tom3}), fits nicely. Ultragraphs are generalizations of graphs, where the range of an edge is a subset of the set of vertices. The associated C*-algebras form a class that strictly contain graph C*-algebras and the Exel-Laca algebras, see \cite{Tom3, TomSimple}. In \cite{MarreroMuhly}, Marrero and Muhly describe ultragraph C*-algebras as groupoid C*-algebras, and the unit space of their ultrapath groupoid was our first choice for the edge shift space associated to an ultragraph. Unfortunately, in its most general case, the description of this topological space is rather technical (a problem also present in the shift space proposed in the work of Exel and Laca), which can make researchers uninterested. 

To avoid the technicalities mentioned above we have used Marrero and Muhly's work as inspiration, and defined ultragraph edge shift spaces in a similar way to the boundary path space defined by Webster. Our space coincide with Webster boundary path space for graphs and also coincide with the edge shift space of a graph when the graph is finite. It also includes new shift spaces, that can not be realized as Webster boundary path space of a graph (see Corollary \ref{newshift}). As expected, our space is intimately related to the path space of the ultragraph. Furthermore, it is always metrizable and possesses a countable basis of clopen sets. For a large class of ultragraphs, that include all graphs, all ultragraphs associated to row finite infinite matrices (and non row finite matrices that satisfy a certain condition), and the ultragraph whose associated C*-algebra is neither a graph algebra nor an Exel-Laca algebra, we show that our space is locally compact and intimately related to ultragraph C*-algebras. In fact, for the class of ultragraphs just mentioned, we realize ultragraph C*-algebras as the partial crossed product of $C_0(X)$ by $\F$, where $\F$ is the free group on the edges of the ultragraph, and $X$ is the shift space we associate to the ultragraph. Building from this realization of ultragraph C*-algebras we show that if two shift spaces are conjugate, via a conjugacy that preserves length, then the associated graph C*-algebras are isomorphic (this shows that our approach to infinite alphabet shift spaces is closely related to ultragraph and graph C*-algebras).

We organize the paper as follows. In section 2 we set up the notation and conventions that will be used through out the paper. Next, in section 3, we define the shift space associated to an ultragraph. We also study the continuity of the shift map $\sigma$, which is not always continuous (see Proposition~\ref{where-shift-cts}), and morphisms between shift spaces. In section 4, we define a partial action on the shift space associated to an ultragraph and show that, for the class of ultragraphs whose shift spaces are locally compact, the associated partial crossed product is isomorphic to the ultragraph C*-algebra. Finally, in section 5, we use the realization of ultragraph C*-algebras as partial crossed products to show that these ultragraph C*-algebras are invariants for conjugacy that preserves length of ultragraph edge shift spaces.

\section{Notation and conventions}

In this section we recall the main definitions and relevant results regarding ultragraphs, as introduced by Tomforde in \cite{Tom3}. We also set up the basic notation we shall use for graphs and ultragraphs, following closely the notation in \cite{MarreroMuhly}. 

\begin{definicao}\label{def of ultragraph}
An \emph{ultragraph} is a quadruple $\mathcal{G}=(G^0, \mathcal{G}^1, r,s)$ consisting of two countable sets $G^0, \mathcal{G}^1$, a map $s:\mathcal{G}^1 \to G^0$, and a map $r:\mathcal{G}^1 \to P(G^0)\setminus \{\emptyset\}$, where $P(G^0)$ stands for the power set of $G^0$.
\end{definicao}

\begin{exemplo}\label{exemplo1} Let $\mathcal{G}$ be the ultragraph with $G^0=\{v_i\}_{i\in \N}$, and edges such that $s(e_i)= v_i$ for all $i$, $r(e_1)=\{v_3, v_4, v_5, \ldots \}$ and $r(e_j)= G^0$ for all $j\neq 1$. We can represent this ultragraph as in the picture below. 
\end{exemplo}

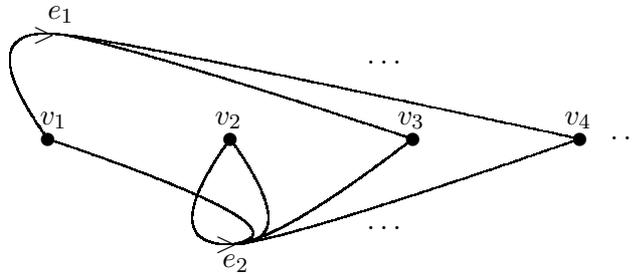
\centerline{
\setlength{\unitlength}{2cm}
\begin{picture}(4,1.2)
\put(0.5,0){\circle*{0.08}}
\put(0.5,0){\qbezier(0,0)(-0.5,0.7)(0,0.7)}
\put(0.4,0.64){$>$}
\put(0.5,0.7){\qbezier(0,0)(0.3,0)(2.4,-0.7)}
\put(0.5,0.7){\qbezier(0,0)(0.3,0)(3.5,-0.7)}
\put(1.7,0){\qbezier(0,0)(-0.5,-0.7)(0,-0.7)}
\put(1.7,0){\qbezier(0,0)(0.5,-0.7)(0,-0.7)}
\put(1.6,-0.75){$>$}
\put(1.7,-0.7){\qbezier(0,0)(0.3,0)(1.2,0.7)}
\put(1.7,-0.7){\qbezier(0,0)(0.3,0)(2.3,0.7)}
\put(1.7,-0.7){\qbezier(0,0)(0.6,0.085)(-1.2,0.7)}
\put(0.45,0.1){$v_{1}$}
\put(1.7,0){\circle*{0.08}}
\put(1.6,0.1){$v_2$}
\put(2.9,0){\circle*{0.08}}
\put(2.8,0.1){$v_3$}
\put(4,0){\circle*{0.08}}
\put(3.9,0.1){$v_4$}
\put(4.2,0){\dots}
\put(2.6,0.5){\dots}
\put(2.6,-0.6){\dots}
\put(0.5,0.8){$e_1$}
\put(1.65,-0.85){$e_2$}
\end{picture}}
\vspace{2 cm}

Before we define the C*-algebra associated to an ultragraph we need the following notion.

\begin{definicao}\label{def of mathcal{G}^0}
Let $\mathcal{G}$ be an ultragraph. Define $\mathcal{G}^0$ to be the smallest subset of $P(G^0)$ that contains $\{v\}$ for all $v\in G^0$, contains $r(e)$ for all $e\in \mathcal{G}^1$, and is closed under finite unions and non-empty finite intersections.
\end{definicao}

\begin{lema}\label{description}\cite[Lemma~2.12]{Tom3}
If $\mathcal{G} ( G^0,
\mathcal{G}^1,r,s)$ is an ultragraph, then \begin{align*} \mathcal{G}^0 = \{
\bigcap_{e
\in X_1} r(e)
\cup \ldots 
\cup \bigcap_{e \in X_n} r(e) \cup F : & \ \text{$X_1,
\ldots, X_n$ are finite subsets of $\mathcal{G}^1$} \\ & \text{ and $F$
is a finite subset of $G^0$} \}.
\end{align*}  
Furthermore, $F$ may be chosen to
be disjoint from $\bigcap_{e \in X_1} r(e) \cup \ldots 
\cup \bigcap_{e \in X_n} r(e)$.
\end{lema}

\begin{definicao}\label{def of C^*(mathcal{G})}
Let $\mathcal{G}$ be an ultragraph. The \emph{ultragraph algebra} $C^*(\mathcal{G})$ is the universal $C^*$-algebra generated by a family of partial isometries with orthogonal ranges $\{s_e:e\in \mathcal{G}^1\}$ and a family of projections $\{p_A:A\in \mathcal{G}^0\}$ satisfying
\begin{enumerate}
\item\label{p_Ap_B=p_{A cap B}}  $p_\emptyset=0,  p_Ap_B=p_{A\cap B},  p_{A\cup B}=p_A+p_B-p_{A\cap B}$, for all $A,B\in \mathcal{G}^0$;
\item\label{s_e^*s_e=p_{r(e)}}$s_e^*s_e=p_{r(e)}$, for all $e\in \mathcal{G}^1$;
\item $s_es_e^*\leq p_{s(e)}$ for all $e\in \mathcal{G}^1$; and
\item\label{CK-condition} $p_v=\sum\limits_{s(e)=v}s_es_e^*$ whenever $0<\vert s^{-1}(v)\vert< \infty$.
\end{enumerate}
\end{definicao}

\subsection{Notation}

Let $\mathcal{G}$ be an ultragraph. A \textit{finite path} in $\mathcal{G}$ is either an element of $\mathcal{G}%
^{0}$ or a sequence of edges $e_{1}\ldots e_{k}$ in $\mathcal{G}^{1}$ where
$s\left(  e_{i+1}\right)  \in r\left(  e_{i}\right)  $ for $1\leq i\leq k$. If
we write $\alpha=e_{1}\ldots e_{k}$, the length $\left|  \alpha\right|  $ of
$\alpha$ is just $k$. The length $|A|$ of a path $A\in\mathcal{G}^{0}$ is
zero. We define $r\left(  \alpha\right)  =r\left(  e_{k}\right)  $ and
$s\left(  \alpha\right)  =s\left(  e_{1}\right)  $. For $A\in\mathcal{G}^{0}$,
we set $r\left(  A\right)  =A=s\left(  A\right)  $. The set of
finite paths in $\mathcal{G}$ is denoted by $\mathcal{G}^{\ast}$. An \textit{infinite path} in $\mathcal{G}$ is an infinite sequence of edges $\gamma=e_{1}e_{2}\ldots$ in $\prod \mathcal{G}^{1}$, where
$s\left(  e_{i+1}\right)  \in r\left(  e_{i}\right)  $ for all $i$. The set of
infinite paths  in $\mathcal{G}$ is denoted by $\mathfrak
{p}^{\infty}$. The length $\left|  \gamma\right|  $ of $\gamma\in\mathfrak
{p}^{\infty}$ is defined to be $\infty$. A vertex $v$ in $\mathcal{G}$ is
called a \emph{sink} if $\left|  s^{-1}\left(  v\right)  \right|  =0$ and is
called an \emph{infinite emitter} if $\left|  s^{-1}\left(  v\right)  \right|
=\infty$. 

For $n\geq1,$ we define
$\mathfrak{p}^{n}:=\{\left(  \alpha,A\right)  :\alpha\in\mathcal{G}^{\ast
},\left\vert \alpha\right\vert =n,$ $A\in\mathcal{G}^{0},A\subseteq r\left(
\alpha\right)  \}$. We specify that $\left(  \alpha,A\right)  =(\beta,B)$ if
and only if $\alpha=\beta$ and $A=B$. We set $\mathfrak{p}^{0}:=\mathcal{G}%
^{0}$ and we let $\mathfrak{p}:=\coprod\limits_{n\geq0}\mathfrak{p}^{n}$. We embed the set of finite paths $\GG^*$ in $\mathfrak{p}$ by sending $\alpha$ to $(\alpha, r(\alpha))$. We
define the length of a pair $\left(  \alpha,A\right)  $, $\left\vert \left(
\alpha,A\right)  \right\vert $, to be the length of $\alpha$, $\left\vert
\alpha\right\vert $. We call $\mathfrak{p}$ the \emph{ultrapath space}
associated with $\mathcal{G}$ and the elements of $\mathfrak{p}$ are called
\emph{ultrapaths}. Each $A\in\mathcal{G}^{0}$ is regarded as an ultrapath of length zero and can be identified with the pair $(A,A)$. We may extend the range map $r$ and the source map $s$ to
$\mathfrak{p}$ by the formulas, $r\left(  \left(  \alpha,A\right)  \right)
=A$, $s\left(  \left(  \alpha,A\right)  \right)  =s\left(  \alpha\right)
$ and $r\left(  A\right)  =s\left(  A\right)  =A$. 

We concatenate elements in $\mathfrak{p}$ in the following way: If $x=(\alpha,A)$ and $y=(\beta,B)$, with $|x|\geq 1, |y|\geq 1$, then $x\cdot y$ is defined if and only if
$s(\beta)\in A$, and in this case, $x\cdot y:=(\alpha\beta,B)$. Also we
specify that:
\begin{equation}
x\cdot y=\left\{
\begin{array}
[c]{ll}%
x\cap y & \text{if }x,y\in\mathcal{G}^{0}\text{ and if }x\cap y\neq\emptyset\\
y & \text{if }x\in\mathcal{G}^{0}\text{, }\left|  y\right|  \geq1\text{, and
if }x\cap s\left(  y\right)  \neq\emptyset\\
x_{y} & \text{if }y\in\mathcal{G}^{0}\text{, }\left|  x\right|  \geq1\text{,
and if }r\left(  x\right)  \cap y\neq\emptyset
\end{array}
\right.  \label{specify}%
\end{equation}
where, if $x=\left(  \alpha,A\right)  $, $\left|  \alpha\right|  \geq1$ and if
$y\in\mathcal{G}^{0}$, the expression $x_{y}$ is defined to be $\left(
\alpha,A\cap y\right)  $. Given $x,y\in\mathfrak{p}$, we say that $x$ has $y$ as an initial segment if
$x=y\cdot x^{\prime}$, for some $x^{\prime}\in\mathfrak{p}$, with $s\left(
x^{\prime}\right)  \cap r\left(  y\right)  \neq\emptyset$. 

We extend the source map $s$ to $\mathfrak
{p}^{\infty}$, by defining $s(\gamma)=s\left(  e_{1}\right)  $, where
$\gamma=e_{1}e_{2}\ldots$. We may concatenate pairs in $\mathfrak{p}$, with
infinite paths in $\mathfrak{p}^{\infty}$ as follows. If $y=\left(
\alpha,A\right)  \in\mathfrak{p}$, and if $\gamma=e_{1}e_{2}\ldots\in
\mathfrak{p}^{\infty}$ are such that $s\left(  \gamma\right)  \in r\left(
y\right)  =A$, then the expression $y\cdot\gamma$ is defined to be
$\alpha\gamma=\alpha e_{1}e_{2}...\in\mathfrak{p}^{\infty}$. If $y=$
$A\in\mathcal{G}^{0}$, we define $y\cdot\gamma=A\cdot\gamma=\gamma$ whenever
$s\left(  \gamma\right)  \in A$. Of course $y\cdot\gamma$ is not defined if
$s\left(  \gamma\right)  \notin r\left(  y\right)  =A$. 

\begin{remark} When no confusion arises we will omit the dot in the notation of concatenation defined above, so that $x\cdot y$ will be denoted by $xy$.
\end{remark}

\begin{definicao}
\label{infinte emitter} For each subset $A$ of $G^{0}$, let
$\varepsilon\left(  A\right)  $ be the set $\{ e\in\mathcal{G}^{1}:s\left(
e\right)  \in A\}$. We shall say that a set $A$ in $\mathcal{G}^{0}$ is an
\emph{infinite emitter} whenever $\varepsilon\left(  A\right)  $ is infinite.
\end{definicao}


\section{Ultragraph shift spaces}

{\bf Throughout assumption:} From now on all ultragraphs in this paper are assumed to have no sinks. Also, $\mathcal{G}$ will always denote an ultragraph.

In this section we will define a shift space associated to a (possibly infinite) ultragraph. In the finite alphabet case, a shift space consists of a set of infinite words with a map (the shift map) that represents the evolution of a discrete system. Since we are dealing with possible infinite alphabets our set will also contain finite sequences. We introduce the precise notions below.

\subsection{The topological space}\label{The topological space}


In this subsection we introduce a topological space associated to an ultragraph. This space generalizes the path space of a directed graph described in \cite{Webster}. Before we define our space we need the following definition.

\begin{definicao}\label{minimal} Let $\GG$ be an ultragraph and $A\in \GG^0$. We say that $A$ is a minimal infinite emitter if it is an infinite emitter that contains no proper subsets (in $\GG^0$) that are infinite emitters. Equivalently, $A$ is a minimal infinite emitter if it is an infinite emitter and has the property that, if $B\in \GG^0$ is an infinite emitter, and $B\subseteq A$, then $B=A$. For a finite path $\alpha$ in $\GG$, we say that $A$ is a minimal infinite emitter in $r(\alpha)$ if $A$ is a minimal infinite emitter and $A\subseteq r(\alpha)$.
We denote the set of all minimal infinite emitters in $r(\alpha)$ by $M_\alpha$.
\end{definicao}

\begin{remark} If $A$ is a minimal infinite emitter and $B$ is an infinite emitter then, since $\GG$ has no sinks, either $A\subseteq B$ or their intersection is at most finite.
\end{remark}


The following result is important in the understanding of minimal infinite emitters.

\begin{lema}\label{miniftyem} Let $x=(\alpha, A) \in \mathfrak{p}$ and suppose that $A$ is a minimal infinite emitter. If the cardinality of $A$ is finite then it is equal to one and, if the cardinality of $A$ is infinite, then $A = \bigcap\limits_{e\in Y} r(e)$ for some finite set $Y\subseteq \GG^1$.
\end{lema}

\demo If the cardinality of $A$ is finite then there exists some $v\in A$ such that $v$ is an infinite emitter. Hence $A=\{v\}$ (since $\{v\}$ is a minimal infinite emitter). 

Suppose that the cardinality of $A$ is infinite. From Lemma \ref{description} we get that $A=\bigcap\limits_{e\in Y_1}r(e)\cup...\cup\bigcap\limits_{e\in Y_n}r(e)\cup F$, where $Y_1,...,Y_n$ are finite subsets of $\GG^1$ and $F\subseteq G^0$ is finite. Note that one of the sets $\bigcap\limits_{e\in Y_i}r(e)$, with $i\in \{1,...,n\}$, or $F$ must be an infinite emitter. Since $A$ is minimal then $A=\bigcap\limits_{e\in Y_i}r(e)$, for some $i\in \{1,...,n\}$, or $A=F$. Since we are assuming that the cardinality of $A$ is infinite, and $F$ is finite, we have that $A\neq F$. Therefore $A=\bigcap\limits_{e\in Y_i}r(e)$ for some $i\in \{1,...,n\}$.

\fim

To construct our space we need to consider the set:
$$X_{fin} = \{(\alpha,A)\in \mathfrak{p}: |\alpha|\geq 1 \text{ and } A\in M_\alpha \}\cup
 \{(A,A)\in \GG^0: A \text{ is a minimal infinite emitter}\}. $$

Let $$X= \mathfrak{p}^{\infty} \cup X_{fin}.$$

To define a basis for a topology in $X$ we need the following notation.

For each $(\beta,B)\in \mathfrak{p}$ let $$D_{(\beta,B)}= \{(\beta, A): A\subseteq B \text{ and } A\in M_\beta \}\cup\{y \in X: y = \beta \gamma', s(\gamma')\in B\}.$$

Suppose $(\beta,B)\in X_{fin}$ and let $F$ be a finite subset of $\varepsilon\left( B \right)$. Define $$D_{(\beta, B),F}=  \{(\beta, B)\}\cup\{y \in X: y = \beta \gamma', \gamma_1' \in \ \varepsilon\left( B \right)\setminus F\}.$$

In the next proposition we describe a basis for the relevant topology in $X$.


\begin{proposicao}\label{basistopology} The collection $\{D_{(\beta,B)}: (\beta,B) \in \mathfrak{p}, |\beta|\geq 1\ \} \cup \{D_{(\beta, B),F}:(\beta, B) \in X_{fin}, F\subseteq \varepsilon\left( B \right), |F|<\infty \}$ 
defined above is a countable basis for a topology on $X$. Furthermore, if $\gamma = e_1 e_2 \ldots \in X$ then a neighborhood basis for $\gamma$ is given by $$\{D_{(e_1 \ldots e_n, r(e_n))}: n\in \N \} $$ and if $x=(\alpha, A) \in X$ then a neighborhood basis for $x$ is given by $$ \{D_{(\alpha, A),F}: F\subseteq \varepsilon\left(  A\right), |F|<\infty \}.$$
\end{proposicao}

\demo Notice that if $\gamma \in \mathfrak{p}^\infty$ is in the intersection, $A\cap B$, of two sets in the collection then it is possible to find an initial segment of $\gamma$, say $\gamma_1$ (so $\gamma = \gamma_1 \gamma^{'}$), such that the cylinder $D_{(\gamma_1,r(\gamma_1))}$ is contained in $A \cap B$. 

Let $x= (\alpha,A) \in X_{fin}$. Suppose that $x\in D_{(\beta_1, B_1),F_1} \cap D_{(\beta_2, B_2),F_2}$ and $|\beta_2|>|\beta_1|$. Then we have that $(\beta_1,B_1)$ is an initial segment of $(\beta_2,B_2)$ and hence $ D_{(\beta_2, B_2),F_2} \subseteq D_{(\beta_1, B_1),F_1}$. If $x\in D_{(\beta, B),F_1} \cap D_{(\beta, B),F_2}$ then $x\in D_{(\beta, B),F_1\cup F_2}$. 
The cases  $x\in D_{(\beta_1, B_1)} \cap D_{(\beta_2, B_2),F_2}$ and  $x\in D_{(\beta_1, B_1)} \cap D_{(\beta_2, B_2)}$ are handled similarly.

For the second part notice that if $\gamma = e_1 e_2 \ldots \in X$ is contained in $D_{(e_1\ldots e_n, B),F}$ for some $n, \ F$ then $e_{n+1}\notin F$ and hence $\gamma \in D_{(e_1\ldots e_{n+1},r(e_{n+1}))}\subseteq D_{(e_1\ldots e_n, B),F}$. If $x=(\alpha, A) \in X$ is contained is some $D_{(\beta, B)}$ then $\beta$ is an initial segment of $\alpha$ and $D_{(\alpha,A)}\subseteq D_{(\beta, B)}$. Finally, if $x$ is contained in some $D_{(\beta, B),F}$, and $\beta \neq \alpha$, then again $\beta$ is an initial segment of $\alpha$ and $D_{(\alpha,A)}\subseteq D_{(\beta, B),F}$. \fim



 
\begin{remark} In the case $\GG$ is a graph then the only possible minimal infinite emitters are sets consisting of a single vertex that, in the graph, is a singular vertex. In this case we can identify an element $(\alpha, r(\alpha))$ with $\alpha$ and one can check that our topological space $X$ coincides with the boundary path space of a graph defined in \cite[Definition~2.1]{Webster}. 
\end{remark}

Next we develop the topological properties of the space $X$ with the topology given by the basis described in Proposition~\ref{basistopology}.

\begin{proposicao}\label{closedbasis} Each basis element given in Proposition \ref{basistopology} is closed.
\end{proposicao}

\demo First we prove that each $D_{(\beta,B)}$, $|\beta|\geq1$, is closed. Suppose $\beta = \beta_1 \ldots \beta_n$. Let $\gamma \in D_{(\beta,B)}^c$.

If $|\gamma| = \infty$, say $\gamma = \gamma_1 \gamma_2 \ldots$ then either $s(\gamma_{n+1}) \notin B$ or $\gamma_1 \ldots \gamma_n \neq \beta$. In both cases $\gamma \in D_{(\gamma_1\ldots \gamma_{n+1}, \{r(\gamma_{n+1})\})} \subseteq D_{(\beta,B)}^c$.

If $|\gamma| < \infty$ and $\gamma = (\gamma_1 \ldots \gamma_k, C)\in X_{fin}$, then we have three cases:

\begin{itemize}
\item $k<n$

We only need to handle the case when $\gamma_1 \ldots \gamma_k$ is an initial subpath of $\beta$. If $s(\beta_{k+1}) \in C$ then $\gamma \in D_{(\gamma_1 \ldots \gamma_k, C),\{\beta_{k+1}\}}\subseteq D_{(\beta,B)}^c$ and if $s(\beta_{k+1}) \notin C$ then $\gamma \in D_{(\gamma_1 \ldots \gamma_k, C)}\subseteq D_{(\beta,B)}^c$.
 
\item $k=n$

The only case we need to deal is when $\gamma_1 \ldots \gamma_k = \beta$, since if $\gamma_1 \ldots \gamma_k \neq \beta$ then $D_{(\gamma_1 \ldots \gamma_k , r(\gamma))}$ is an open neighborhood of $\gamma$. So, suppose that $\gamma_1 \ldots \gamma_k = \beta$. We have two cases. 

If $|C|< \infty$ then, by Lemma \ref{miniftyem}, $|C|=1$ and we can write $C=\{ v \}$. This implies that $v\notin B$ since otherwise we would have $\gamma \in D_{(\beta, B)}$. Hence $D_\gamma \subseteq  D_{(\beta,B)}^c$. 

Now suppose that $|C| = \infty$. Then $|C\cap B|< \infty$ is finite, since otherwise the minimality of $C$ implies that $C\cap B = C$ and hence $C\subseteq B$ and $\gamma \in D_{(\beta,B)}$, a contradiction. Notice also that the vertices in $B\cap C$ are not infinite emitters, that is,  if $v \in B\cap C$ then $|s^{-1}(v)|<\infty$, otherwise $C$ is not a minimal infinite emitter. Let $F= \{e \in \GG^1: s(e) \in B\cap C\}$. Then $F$ is finite and $D_{\gamma, F} \subseteq D_{(\beta,B)}^c$.

\item $k>n$

Clearly the only case we need to deal is when $\gamma_1 \ldots \gamma_n = \beta_1 \ldots \beta_n$. In this case $s(\gamma_{n+1})\notin B$ and hence $\gamma\in D_{(\gamma_1\ldots \gamma_{n+1},r(\gamma_{n+1}))}\subseteq D_{(\beta,B)}^c$.
\end{itemize}


Next we prove that each $D_{(B, B)}$, with $B$ a minimal infinite emitter, is closed. 

Let $x \in D_{(B,B)}^c$. Suppose that $|x|=0$, that is, $x=(A,A)$. If $A\cap B = \emptyset$ then $D_x$ is an appropriate open neighboorhood. If $A \cap B \neq \emptyset$ then $|A\cap B| <\infty$ and, as before, each $v\in A \cap B$ is not an infinite emitter. Let $F= \{e \in \GG^1: s(e) \in A \cap B\}$. Then $D_{(A,A),F} \subseteq D_{(B,B)}^c$. It is straightforward to deal with the case $|x|\geq 1$ .

To finish notice that for each $(\beta,B)\in X_{fin}$ and finite $F \subseteq \ \varepsilon\left( B \right)$ we can write $$D_{(\beta, B),F}= \displaystyle D_{(\beta, B)} \bigcap_{\lambda \in F} D_{(\beta \lambda, r(\lambda))}^c.$$ \fim


\begin{proposicao} $X$ is Hausdorff.
\end{proposicao}

\demo Let $x\neq y $ be two elements of $X$. If $x$ and $y$ are both in $\mathfrak{p}^\infty$, say $x= \gamma_1 \gamma_2 \ldots$ and $y= \beta_1 \beta_2 \ldots $, then there exists $n$ such that $\gamma_n \neq \beta_n$ and hence we can choose $D_{(\beta_1\ldots \beta_n, r(\beta_n))}$ as separating open neighborhoods. 

Suppose now that $x \in \mathfrak{p}^\infty$ and $y \in X_{fin}$, say $x= \gamma_1 \gamma_2 \ldots$ and $y= (\beta, B)$. It is clear how to find separating neighborhoods if $\beta$ is not an initial segment of $x$. So suppose that $\beta$ is an initial segment of $x$, say $\beta = \gamma_1 \ldots \gamma_k$, and $\gamma_{k+1} \in B$ (it is clear how to separate $x$ and $y$ if $\gamma_{k+1} \notin B$). Notice that in this case $B$ is an infinite minimal emitter, and $D_{(\gamma_1 \ldots \gamma_{k+1}, r(\gamma_{k+1}))}$ and $D_{(\beta, B),\{\gamma_{k+1}\}}$ are separating open neighborhoods.

We are left with the case when $x,y \in X_{fin}$, say $x=(\alpha, A)$ and $y=(\beta, B)$. As before there are a few cases to consider. We will show how to obtain separating neighborhoods when $\alpha = \beta$ as other cases are straightforward. If $\alpha= \beta$ and $A\cap B = \emptyset $ then $D_{(\alpha,A)}$ and $D_{(\beta,B)}$ are separating neighborhoods. If $A\cap B \neq \emptyset$ then $|A \cap B|<\infty$ and, since $A$ is a minimal infinite emitter, $ _{A\cap B}\GG^1$ is a finite set. Hence $D_{(\alpha, A), _{A\cap B}\GG^1}$ and  $D_{(\beta, B), _{A\cap B}\GG^1}$ are separating open neighborhoods. \fim

\begin{remark} Since our space is Hausdorff and has a countable basis of clopen sets it follows from Urysohn's metrization Theorem that $X$ is metrizable.
\end{remark}

For metrizable spaces, it is useful to have a description of convergence of sequences. In light of Proposition \ref{basistopology} we have:

\begin{corolario}\label{convseq} Let $\{x^n\}_{n=1}^{\infty}$ be a sequence of elements in $X$, where $x^n = (\gamma^n_1\ldots \gamma^n_{k_n}, A_n)$ or $x^n = \gamma_1^n \gamma_2^n \ldots$, and let $x \in X$.
\begin{enumerate}[(a)]
\item If $|x|= \infty$, say $x=\gamma_1 \gamma_2 \ldots$, then $\{x^n\}_{n=1}^{\infty}$ converges to $x$ if, and only if, for every $M\in \N$ there exists $N\in \N$ such that $n>N$ implies that $|x^n|\geq M$ and $\gamma^n_i= \gamma_i$ for all $1\leq i \leq M$.
\item If $|x|< \infty$, say $x=(\gamma_1 \ldots \gamma_k, A)$, then $\{x^n\}_{n=1}^{\infty}$ converges to $x$ if, and only if, for every finite subset $F\subseteq \varepsilon\left(  A\right)$ there exists $N\in \N$ such that $n > N$ implies that $x^n = x$ or $|x^n|> |x|$, $\gamma^n_{|x|+1} \in \ \varepsilon\left(  A\right)\setminus F$, and $\gamma^n_i = \gamma_i$ for all $1 \leq i \leq |x|$. 
\end{enumerate}
\end{corolario}

It would be convenient if our space has a basis of open, compact sets. But this is not true in general, as we remark below. So we will need to add an extra hypothesis.

\begin{remark} Notice that in general it is not true that every $D_{(\beta, B)}$, $(\beta, B)\in \mathfrak{p}$, is compact. For example, if $\mathcal{G}$ is an ultragraph with an edge $e$ such that $r(e)$ contains an infinite number of vertices, and each of these vertices is an infinite emitter, then $D_{(e,r(e))}$ is not compact.  
\end{remark}

{\bf Condition (RFUM):} We say that an ultragraph $\mathcal{G}$ satisfies condition (RFUM) if for each edge $e\in \GG^1$ its range can be written as $$r(e) = \displaystyle \bigcup_{n=1}^k A_n,$$ where $A_n$ is either a minimal infinite emitter or a single vertex. 

The hypothesis above includes any graph and the ultragraph whose algebra is neither an Exel-Laca nor a graph algebra (see \cite{TomSimple}), but it does not include all ultragraphs associated to infinite matrices. Recall (as in \cite{Tom3}) that if $A$ is an infinite matrix of zeros and ones then the associated Exel-Laca algebra (as in \cite{ExelLaca}) is isomorphic to the ultragraph C*-algebra $C^*(\GG_A)$, where $\GG_A$ is the ultragraph given by $G^0=\{i:i\in\N\}$, $\GG_A^1=\{e_{i}:i\in \N \}$, $s(e_i)=i$ and $r(e_i)=\{j: A_{ij} =1\}$. Any row finite matrix is an example of a matrix such that $\GG_A$ satisfies Condition (RFUM). Other examples (and counter-examples) are:

\begin{exemplo} Let $A$ be the infinite matrix with all entries equal to one, $B= \left( \begin{smallmatrix}
  1 & 1 & 1 & 1 & \cdots \\
  1 & 0 & 1 & 0 & \cdots \\
  0 & 1 & 0 & 1 & \cdots \\
  1 & 0 & 1 & 0 & \cdots \\
  \vdots  & \vdots  & \vdots & \vdots & \ddots 
 \end{smallmatrix} \right),$ whereas $C= \left( \begin{smallmatrix}
  1 & 1 & 1 & 1 & \cdots \\
  1 & 0 & 1 & 0 & \cdots \\
  1 & 0 & 1 & 0 & \cdots \\
  \vdots  & \vdots  & \vdots & \vdots & \ddots 
 \end{smallmatrix} \right).$ 
Then the ultragraphs $\GG_A$ and $\GG_B$ satisfy Condition (RFUM) and $\GG_C$ does not satisfy it.
 \end{exemplo}



\begin{proposicao}\label{compactbasis} Let $\GG$ be an ultragraph that satisfies Condition (RFUM). Then each basis element of the topology on $X$, as in Proposition \ref{basistopology}, is compact.
\end{proposicao}

\demo For each $(\beta, B) \in \mathfrak{p}$, with $|\beta|\geq 1$,
we will show that $D_{(\beta,B)}$ is sequentially compact. Let $(x^n)$ be a sequence in $D_{(\beta,B)}$. Then $x^n = (\beta, A_n)$, $x^n = (\beta \alpha_1^n \ldots \alpha_{k_n}^n, A_n)\in X_{fin}$ or $x^n = \beta \alpha_1^n \alpha_2^n \ldots \in \mathfrak{p}^\infty$. 

By Condition (RFUM) there is only a finite number of minimal infinite emitters contained in $B$. So, if there exists an infinite number of indices such that $x^n = (\beta, A_n)$ we obtain a constant subsequence of $(x^n)$. Therefore we can assume, without loss of generality, that $|x^n|> |\beta|$ for all $n$.

If there is not an infinite number of indices such that $\alpha_1^n$ coincide then we can assume without loss of generality that $\alpha_1^n \neq \alpha_1^m$ for all $n\neq m$. In this case $B$ must contain at least one infinite emitter and hence can be written as a finite union of minimal infinite emitters and unitary sets (this follows from Condition (RFUM)).
 So, we can find $A\in M_\beta$ such that $|\{n : s(\alpha_1^n)\in A \}| = \infty$ and hence we obtain a subsequence converging to $(\beta,A)$.

If there exists an infinite number of indices such that $\alpha_1^n$ coincide, say $\alpha_1^n = \gamma_1$, then we pass to a subsequence such that $x^n = (\beta \gamma_1 \alpha_2^n \ldots \alpha_{k_n}^n, A_n)$ or $x^n = \beta \gamma_1 \alpha_2^n \alpha_3^n \ldots$ and notice that each $x^n \in D_{(\beta \gamma_1, r(\gamma_1))}$. We now repeat the procedure described in the above paragraph and either obtain a subsequence converging to a finite sequence of pass to a subsequence such that $x^n = (\beta \gamma_1 \gamma_2 \alpha_3^n \ldots \alpha_{k_n}^n, A_n)$ or $x^n = \beta \gamma_1 \gamma_2 \alpha_3^n \alpha_4^n \ldots$. Proceeding inductively we either obtain at step k a subsequence converging to a finite sequence or, through a Cantor diagonal argument, obtain a subsequence converging to the infinite sequence $\beta\gamma_1 \gamma_2 \ldots$.

We conclude that $D_{(\beta, B)}$ is sequentially compact. Finally, notice that the proof that the sets of the form $D_{(\beta, B),F}$ are sequentially compact is completely analogous. \fim
 
It is interesting to notice that the set of infinite sequences is dense in $X$, as we show in the next proposition.

\begin{proposicao} The set of infinite sequences, $\mathfrak{p}^\infty$, is dense in $X$.
\end{proposicao}

\demo Let $x=(\alpha,A) \in X_{fin}$. Since $\GG$ has no sinks, for each edge $e_n\in \ \varepsilon\left(  A\right)$ there exists an infinite path $\gamma^n=\gamma^n_1 \gamma^n_2 \ldots$ such that $\gamma^n_1 = e_n$. Hence, by Corollary \ref{convseq}, $x$ is the limit of the sequence $\{\alpha \cdot \gamma^n \}_{n=1}^{\infty}$.\fim
 
\subsection{The shift map}

In this subsection we define the shift map on $X$ and study its continuity.

\begin{definicao}\label{shift-map-def}
The \emph{shift map} is the function $\sigma : X \rightarrow X$ defined by $$\sigma(x) =  \begin{cases} \gamma_2 \gamma_3 \ldots & \text{ if $x = \gamma_1 \gamma_2 \ldots \in \mathfrak{p}^\infty$} \\ (\gamma_2 \ldots \gamma_n,A) & \text{ if $x = (\gamma_1 \ldots \gamma_n,A) \in X_{fin}$ and $|x|> 1$} \\(A,A) & \text{ if $x = (\gamma_1,A) \in X_{fin}$} \\ (A,A) & \text{ if $x = (A,A)\in X_{fin}$.} 
\end{cases}$$
\end{definicao}

Notice that if $|x| = \infty$ then $|\sigma(x)| = \infty$, if $|x| \in \N$ then $|\sigma(x)| = |x|-1$, and if $|x| = 0$ then $|\sigma(x)| = 0$.  

\begin{remark}\label{shiftdefzero} In graph C*-algebra theory it is usual to not define the shift map on elements of length zero (see for example \cite{NTW}), but in symbolic dynamics this is often the case (see \cite{GSS, GSS1, OTW}). We therefore have chosen to define $\sigma$ in all elements of $X$.
\end{remark}

\begin{proposicao}\label{where-shift-cts}
The shift map $\sigma : X \rightarrow X$ is continuous at all points of $X$ with length greater than zero. In addition, if $|x|\geq 1$ then there exists an open set $U$ that contains no elements of length zero such that $x \in U$, $\sigma(U)$ is an open subset of $X$, and $\sigma|_U : U \to \sigma(U)$ is a homeomorphism.
\end{proposicao}

\demo Suppose $\{x^n\}_{n=1}^{\infty}$ is a sequence in $X$, where $x^n = (\gamma^n_1\ldots \gamma^n_{k_n}, A_n)$ or $x^n = \gamma_1^n \gamma_2^n \ldots$, and suppose that $\{x^n\}$ converges to $x \in X$.

If $|x|=\infty$ then it is clear that $\{\sigma(x^n)\}$ converges to $\sigma (x)$ (the proof for the usual shift map on finite alphabets works in this case). So we will focus in the case when $|x|<\infty$. 

Suppose $\{x^n\}$ converges to $x=(\gamma_1,\ldots,\gamma_k,A)$. By Corollary \ref{convseq}, given $F\subseteq \ \varepsilon\left(  A\right)$ there exists $N\in \N$ such that for all $n>N$, $|x^n|> k$, $\gamma^n_j= \gamma_j$ for $1 \leq j \leq k$, and $\gamma^n_{k+1} \in \varepsilon\left(  A\right)\setminus F$, or $x^n = x$. Hence, for all $n>N$, $|\sigma(x^n)|> k-1$ and $\sigma(x^n)_{k} = \gamma^n_{k+1}\in \varepsilon\left(  A\right)\setminus F$, or $\sigma(x^n) = \sigma(x)$. We conclude that $\{\sigma(x^n)\}$ converges to $\sigma(x)$ as desired.

For the second part of the proposition, notice that if $|x|\geq 1$ and $U$ is one of the basis neighborhoods that contain x (given in Proposition~\ref{basistopology}) then $\sigma$ is a homeomorphism between $U$ and $\sigma(U)$. \fim

\subsection{The shift space and its morphisms}

Given an ultragraph $\GG$ we now define the associated shift space.

\begin{definicao} If $\GG$ is an ultragraph, we defined the associated one-sided shift space to be the pair $(X,\sigma)$, where $X$ is the topological space from Section \ref{The topological space} and $\sigma:X\rightarrow X$ is the map from Definition \ref{shift-map-def}. We will often refer to the space $X$ with the understanding that the map $\sigma$ is attached to it.
\end{definicao}

\begin{remark} Notice that if $\GG$ is a finite graph then there are no infinite emitters and $X$ (as a topological space) is the usual infinite path space of the graph. So our definition coincides with the usual definition of a one sided edge shift. In symbolic dynamics it is also usual to consider the associated two-sided shift (which can be seen as an inverse limit of the one sided shift), but in the infinite alphabet case it is not clear which topology to consider in the two-sided shift. Work on two-sided shift spaces over infinite alphabets, related to Ott-Tomforde-Willis one sided shift spaces, was done in \cite{GSS1}.
\end{remark}

Next we define the appropriate morphisms between shift spaces associated to ultragraphs.

\begin{definicao} Let $\GG_1$ and $\GG_2$ be ultragraphs and denote the associated shift spaces by $X_{\GG_1}$ and $X_{\GG_2}$ respectively. We say that $\phi:X_{\GG_1}\rightarrow X_{\GG_2}$ is a shift morphism if it is continuous and shift commuting (that is $\sigma \circ \phi = \phi \circ \sigma$). We say that $\phi$ is length preserving if $|\phi(x)|=|x|$ for all $x\in X_{\GG_1}$. If a shift morphism is also a homeomorphism then we call it a conjugacy and say that $X_{\GG_1}$ and $X_{\GG_2}$ are conjugate.
\end{definicao}

\begin{remark} If we had chosen to define the shift map only on $X\setminus \mathfrak{p}^0$ (see Remark~\ref{shiftdefzero}) then  we could define a shift morphism to be a continuous map $\phi:X_{\GG_1}\rightarrow X_{\GG_2}$ such that $\phi(X_{\GG_1}\setminus \mathfrak{p}^0) \subseteq X_{\GG_2}\setminus \mathfrak{p}^0$ and $\phi$ is shift commuting. Then a conjugacy (a bijection such that both $\phi$ and $\phi^{-1}$ are shift morphisms) would automatically be length preserving.
\end{remark}

We remarked above that our construction of edge shift spaces associated to ultragraphs generalizes the notion of edge shift spaces of graphs. Below we show that the edge shift associated to a finite ultragraph is conjugate to the edge shift of a graph.

\begin{proposicao} 
Let $\mathcal{G}=(G^0, \mathcal{G}^1, r,s)$ be a finite ultragraph. Then there exists a graph $\mathcal{F}$ such that $X_\GG$ and $X_{\mathcal{F}}$ are conjugate (by a length preserving conjugacy).
\end{proposicao}

\demo Let $\mathcal{G}=(G^0, \mathcal{G}^1, r,s)$ be a finite ultragraph. Consider the graph $\mathcal{F}=(G^0, \mathcal{F}^1, r,s)$, where the edges in $\mathcal{F}^1$ are defined, for each $e\in \GG^1$ and $v\in r(e)$, by $s(f_{e_v})= s(e)$ and $r(f_{e_v})= v$. Then $\mathcal{F}^1=\{f_{e_v}: e \in \GG^1, v\in r(e)\}$. We will show that $X_\GG$ and $X_{\mathcal{F}}$ are conjugate. 

Notice that since the $\GG$ is finite $\mathcal{F}$ is also finite. So $X_\GG$ and $X_{\mathcal{F}}$ contain only infinite sequences and our definition of conjugacy becomes the usual definition of conjugacy between shift spaces. 

Define $\phi:X_\GG \rightarrow X_{\mathcal{F}}$ by $\phi(e_1 e_2 e_3 \ldots) = g_1 g_2 g_3 \ldots$, where $g_i$ is the edge $f_{{e_i}_{s(e_{i+1})}}$ , that is, $g_i$ is the graph edge associated to $e_i$ such that $s(g_i) = s(e_i)$ and $r(g_i) = s(e_{i+1})$. Clearly $\phi$ is shift commuting, and since the topology in both shift spaces is given by the usual topology of cylinder sets if follows that $\phi$ is a conjugacy as desired. \fim

\section{Ultragraph C*-algebras as partial crossed products}

In this section we define a partial action on the shift space $X$ associated to an ultragraph with no sinks that satisfy Condition (RFUM), and show that the associated crossed product is isomorphic to the ultragraph C*-algebra. This generalizes a previous result by the authors, see \cite[Theorem 4.11]{GRUltra}, to a larger class of ultragraphs. 
We emphasize the assumption that in this section, and the next one, all ultragraphs satisfy Conditon (RFUM):

{\bf Assumption:} From now on we assume that all ultragraphs satisfy Condition (RFUM).

The results in this section are crucial in our proof that ultragraph C*-algebras are invariants for shift conjugacy. We start with the following:

\begin{lema}\label{cvclopen} For each $v\in G^0$ the set $$X_v=\{(\alpha,A)\in X_{fin}:s(\alpha)=v\}\cup \{\gamma\in \mathfrak{p}^\infty:s(\gamma)=v\}= \{x\in X:s(x)=v\}$$ is nonempty, clopen, and compact. 
\end{lema}

\demo Let $v\in G^0$. Since the ultragraphs we are considering have no sinks we have that $X_v\neq \emptyset$. 

If $v$ is an infinite emitter then $X_v=D_{(\{v\},\{v\})}$ and so $X_v$ is clopen by Proposition \ref{closedbasis}, and compact by Proposition \ref{compactbasis}. If $v$ is a finite emitter then $s^{-1}(v)$ is finite and $X_v=\bigcup\limits_{e\in s^{-1}(v)}D_{(e,r(e))}$. Since, by Propositions \ref{closedbasis} and \ref{compactbasis}, each $D_{(e,r(e))}$ is clopen and compact, it follows that $X_v$ is clopen and compact. \fim

Let $\F$ be the free group generated by $\mathcal{G}^1$. We will define a partial action of $\F$ on $X$. For this, let $P\subseteq \F$ be defined by $$P:=\{e_1...e_n\in \F: e_i\in \mathcal{G}^1: n\geq 1\},$$ 
and definte the sets $X_c$, for each $c\in \F$, as follows:
\begin{itemize}
\item for the neutral element $0\in \F$ let $X_0=X$;

\item for $a\in P$ define \nl$X_a=\{(\beta,B)\in X_{fin}:\beta_1...\beta_{|a|}=a\}\cup\{\gamma\in \mathfrak{p}^\infty: \gamma_1...\gamma_{|a|}=a\}$; and

$X_{a^{-1}}=\{(A,A)\in X_{fin} :A\subseteq r(a)\}\cup\nl\cup\{(\beta,B)\in X_{fin}: s(\beta)\in r(a)\}\cup \{\gamma\in \mathfrak{p}^\infty:s(\gamma)\in r(a)\};$ 

\item for $a,b\in P$ with $ab^{-1}$ in its reduced form, define \nl $X_{ab^{-1}}=\left\{(a,A)\in X_{fin}:A\subseteq r(a)\cap r(b)\right\}\cup\nl\cup\left\{(\beta,B)\in X_{fin}:\beta_1...\beta_{|a|}=a \text{ and } s(\beta_{|a|+1})\in r(a)\cap r(b)\right\}\cup\nl\cup\left\{\gamma \in \mathfrak{p}^\infty:\gamma_1...\gamma_{|a|}=a \text{ and } s(\gamma_{|a|+1})\in r(a)\cap r(b)\right\}$;

\item for all other $c\in \F$ define $X_c=\emptyset$.

\end{itemize}

\begin{remark} Notice that if $a\in P$ is not a path then $X_a$ is empty. Analogously, if $a,b\in P$ and $r(a)\cap r(b) = \emptyset$ then $X_{ab^{-1}}$ is empty.
\end{remark}

To obtain a topological partial action we need to prove that each $X_c$ is open. In fact, we have

\begin{proposicao}\label{opencompact}
The subsets $X_c$, with $c\in \F$, are clopen and compact.
\end{proposicao}

\demo Let $a\in P$. Notice that $X_a=D_{(a,r(a))}$ and hence it is open. Also, $X_a$ is closed by Proposition \ref{closedbasis}, and compact by Proposition \ref{compactbasis}. Now we turn to $X_{a^{-1}}$. Notice that, by the Condition (RFUM), $r(a)=r(a_{|a|})=\bigcup\limits_{n=1}^k A_n$, where each $A_n$ in a minimal infinite emitter or a single vertex. If $A_n$ is a minimal infinite emitter then $D_{(A_n,A_n)}$ is clopen by Proposition \ref{closedbasis} and compact by Proposition \ref{compactbasis}. If $A_n$ is a single vertex then $D_{(A_n,A_n)}$ is clopen and compact by Lemma \ref{cvclopen}. Since $X_{a^{-1}}=\bigcup\limits_{n=1}^k D_{(A_n,A_n)}$ we obtain that $X_{a^{-1}}$ is clopen and compact.

Next we prove that, for each $a,b\in P$, the set $X_{ab^{-1}}$ is closed and compact. Let $(\beta^n)_{n\in \N}\subseteq X_{ab^{-1}}$ be a sequence converging to $\beta$. Since $X_{ab^{-1}}\subseteq X_a$, and $X_a$ is closed, we have that $\beta \in X_a$. We now divide the proof on two steps, depending on the length of $\beta$:

Suppose that $|\beta|>|a|$. Let $n_0$ be such that $\beta^n \in D_{(\beta_1...\beta_{|a|+1}, r(\beta_{|a|+1}))}$ for all $n\geq n_0$. Then $a=\beta_1^{n_0}...\beta_{|a|}^{n_0}=\beta_1...\beta_{|a|}$ and $\beta_{|a|+1}^{n_0}=\beta_{|a|+1}$. Since $s(\beta_{|a|+1}^{n_0})\in r(a)\cap r(b)$, we have that $\beta \in X_{ab^{-1}}$. 

Now, suppose that $|\beta|=|a|$. Then $\beta=(a,A)$, with $A\subseteq r(a)$ a minimal infinite emitter. We need to show that $A\subseteq r(b)$ (so that $(a,A)\in X_{ab^{-1}}$). By Proposition \ref{convseq} there exists $n_0$ such that, for all $n\geq n_0$, either $\beta^n=\beta$ or $|\beta^n|>|a|$ with $\beta_1^{n}...\beta_{|a|}^n=a$ and $\beta_{|a|+1}^n\in \varepsilon\left(  A\right)$.
If $\beta^n=\beta$ for some $n$ then we are done and so we can assume, without loss of generality, that $|\beta^n|>|a|$ for all $n\geq n_0$. Suppose that $\{\beta_{|a|+1}^n:n\geq n_0\}$ is finite, and denote this set by $F$. Then $\beta^n\notin D_{(a,A),F}$ for all $n\geq n_0$, which by Proposition $\ref{convseq}$ is not possible, since $\beta^n\rightarrow (a,A)$. So the set $\{\beta_{|a|+1}^n:n\geq n_0\}$ is infinite. Since $\{s(\beta_{|a|+1}^n):n\geq n_0\}\subseteq r(b)\cap r(a)$ and $\{s(\beta_{|a|+1}^n):n\geq n_0\}\subseteq A$ then $\{s(\beta_{|a|+1}^n):n\geq n_0\}\subseteq A\cap r(b)$ and hence $A\cap r(b)$ is an infinite emitter. Since $A$ is minimal we have that $A=A\cap r(b)$. So, it follows that $A\subseteq r(b)$ and hence $\beta=(a,A)\in X_{ab^{-1}}$.

We conclude that $X_{ab^{-1}}$ is closed. Since $X_a$ is compact and $X_{ab^{-1}}\subseteq X_a$ it follows that $X_{ab^{-1}}$ is compact.

Finally we show that $X_{ab^{-1}}$ is open. Let $\beta\in X_{ab^{-1}}$. If $\beta \in \mathfrak{p}^\infty$, then $\beta=a\beta'$ with $s(\beta')\in r(a)$, and $\beta\in D_{(a\beta_1', r(\beta_1'))}\subseteq X_{ab^{-1}}$. If $\beta=(\alpha,A)$ with $|\alpha|>|a|$, then $\alpha=a\alpha'$, where $s(\alpha')\in r(a)$, and hence $\beta\in D_{(a\alpha', r(\alpha'))}\subseteq X_{ab^{-1}}$. If $|\beta|=|a|$ then $\beta=(a,A)$, where $A\subseteq r(a)\cap r(b)$ is a minimal infinite emitter. Note that $\beta\in D_{(a,A)}\subseteq X_{ab^{-1}}$ and hence each element of $X_{ab^{-1}}$ is an interior element. Therefore $X_{ab^{-1}}$ is open. \fim 

The next step to define a topological partial action is to define maps between the non-empty sets $X_c$, $c\in \F$. We do this below.

For each $a\in P$ such that $X_a$ is non-empty, let $\theta_a:X_{a^{-1}}\rightarrow X_a$ be defined by $\theta_a(x)=a \cdot x$, for each $x\in X_{a^{-1}}$ (here we are using the embedding of $a$ in $\mathfrak{p}$ as $(a,r(a))$). More explicitly, $\theta_a(x)$ is defined in the following way: $\theta_a((A,A))=(a,A)$, $\theta_a(\beta,B)=(a\beta,B)$ and $\theta_a(\gamma)=a\gamma$.
Moreover, define $\theta_{a^{-1}}:X_a\rightarrow X_{a^{-1}}$ by $\theta_{a^{-1}}(y)=\widehat{a}y$, where we understand $\widehat{a}y$ as
$\theta_{a^{-1}}((a,A))=(A,A)$, $\theta_{a^{-1}}(ab,B)=(b,B)$ and $\theta(a\gamma)=\gamma$. Notice that $\theta_a$ and $\theta_{a^{-1}}$ are inverses each of other. Finally, for $a,b\in P$ such that $X_{ab^{-1}}$ is non-empty  define $\theta_{ab^{-1}}:X_{ba^{-1}}\rightarrow X_{ab^{-1}}$ by $\theta_{ab^{-1}}(x)=a\widehat{b}x$, so that 
$\theta_{ab^{-1}}(b,A)=(a,A)$, $\theta_{ab^{-1}}(b\alpha,B)=(a\alpha,B)$ and $\theta_{ab^{-1}}(b\gamma)=a\gamma$.

\begin{proposicao} For each $c\in \F$, the map $\theta_c:X_{c^{-1}}\rightarrow X_c$, as defined above, is a homeomorphism.
\end{proposicao}

\demo Let $a,b\in P$. We show that $\theta_{ab^{-1}}:X_{ba^{-1}}\rightarrow X_{ab^{-1}}$ is continuous. Let $(x^n)_{n\in \N}\subseteq X_{ba^{-1}}$ be such that $x^n\rightarrow x\in X_{ba^{-1}}$. Note that $x^n=b\gamma_1^n\gamma_2^n...$ or $x^n=(b\gamma_1^n...\gamma_{k_n}^n,A_n)$ for each $n\in \N$. If $x\in \mathfrak{p}^\infty$, say $x=b\gamma_1\gamma_2...$,  then for each $M\in \N$ there exists $n_0\in \N$ such that $|x^n|>M+|b|$ and $b\gamma_1^n...\gamma_M^n=b\gamma_1\gamma_2...\gamma_{M}$, for each $n\geq n_0$. Let $y^n=\theta_{ab^{-1}}(x^n)$ for each $n$ and $y=\theta_{ab^{-1}}(x)$. Then $y^n=(a\gamma_1^n...\gamma_{k_n}^n,A_n)$ or $y^n=a\gamma_1^n\gamma_2^n...$ for each $n$. Since $a\gamma_1^n...\gamma_M^n=a\gamma_1...\gamma_M$, for each $n\geq n_0$ then by Corollary \ref{convseq}, $y^n\rightarrow a\gamma_1\gamma_2...=y$. If $|x|<\infty$ then (using Corollary \ref{convseq} again) we obtain that $\theta_{ab^{-1}}(x^n)\rightarrow \theta_{ab^{-1}}(x)$.

The inverse of $\theta_{ab^{-1}}$ is $\theta_{ba^{-1}}$ which is also continuous, and so $\theta_{ab^{-1}}$ is a homeomorphism.

The verification that $\theta_a:X_{a^{-1}}\rightarrow X_a$ is a homeomorphism, for each $a\in P$, follows in a similar way to what was done in the previous case. \fim

\begin{corolario} Since, for each $t\in \F$, the map $\theta_t:X_{t^{-1}}\rightarrow X_t$ is a homeomorphism we get that $\alpha_t:C(X_{t^{-1}})\rightarrow C(X_t)$, defined by $\alpha_t(f)=f\circ \theta_{t^{-1}}$, is a *-isomorphism. It follows from the definitions of $X_g$ and $\theta_g$, for $g\in \F$, that $\theta_t(X_{t^{-1}}\cap X_h)\subseteq X_{th}$ and that $\theta_t(\theta_h(x))=\theta_{th}(x)$ for each $x\in X_{(th)^{-1}}\cap X_{h^{-1}}$. Hence $(\{X_t,\theta_t)\}_{t\in \F}$ is a topological partial action. Consequently, $(\{C(X_t)\}_{t\in \F}, \{\alpha_t\}_{t\in \F})$ is a C*-algebraic partial action of $\F$ in $C(X)$ (see for example \cite{beuter, exelbook}). Note that if $t\in \F$ is not of the form $t=ab^{-1}$, with $a,b\in P\cup\{0_{\F}\}$ then $C(X_t)$ is the null ideal and $\alpha_t$ is the null map.
\end{corolario}

\begin{remark} The partial action above coincides with the one introduced in \cite{toke} when $\mathcal{G}$ is a directed graph without sinks.
\end{remark}

To prove that the C*-algebra of an ultragraph $\GG$, $C^*(\GG)$, is isomorphic to the partial crossed product associated to the partial action defined above we will need to specify the image, by the isomorphism, of elements in $C^*(\GG)$ of the form $P_A$, where $A\in \GG^0$. For this we need the definition below and a couple of lemmas.

\begin{definicao} For each $A\in \mathcal{G}^0$ define $X_A\subseteq X$ by $$X_A=\{x\in X:s(x)\subseteq A\}.$$
\end{definicao}

\begin{lema}\label{lemaAB} For each $A,B\in \mathcal{G}^0$ it holds that $X_{A\cap B}=X_A\cap X_B$ and $X_{A\cup B}=X_A\cup X_B$.
\end{lema}

\demo Let $A,B\in \mathcal{G}^0$. The equality $X_{A\cap B}=X_A\cap X_B$ follows directly from the definition of the sets. We show that $X_{A\cup B}=X_A\cup X_B$. It is clear that $X_{A}\cup X_{B}\subseteq X_{A\cup B}$, since $A\subseteq A\cup B$ and $B\subseteq A\cup B$. We prove the other containment. Let $x\in X_{A\cup B}$, that is, $s(x)\subseteq A\cup B$. If $|x|\geq 1$ then $x=(\beta, B)$ and $s(x)=s(\beta)$, which is a single vertex, and so $s(x)\in A$ or $s(x)\in B$, that is, $x\in X_A\cup X_B$. If $|x|=0$ then $x=(D,D)$, where $D\in \mathcal{G}^0$ is a minimal infinite emitter. Since $D\subseteq A\cup B$ then $D=(D\cap A)\sqcup (D\cap (B\setminus A))$, and hence $D\cap A$ or $D\cap (B\setminus A)$ is an infinite emitter. Suppose (without loss of generality) that $D\cap A$ is an infinite emitter. Note that $D\cap (B\setminus A)=\emptyset$, since otherwise $D\cap A\subsetneqq D$, which is impossible since $D$ is a minimal infinite emitter. It follows that $D=D\cap A$, that is, $D\subseteq A$, and so $x=(D,D)\in X_A$. Therefore $X_{A\cup B}\subseteq X_A\cup X_B$. \fim

\begin{corolario} Each $X_A$ is clopen and compact.
\end{corolario}

\demo First note that $X_{r(e)}=X_{e^{-1}}$ for each $e\in \mathcal{G}^1$, and by Proposition \ref{opencompact}, $X_{e^{-1}}$ is clopen and compact. Moreover, for each $v\in G^0$, $X_v$ is clopen and compact by Lemma $\ref{cvclopen}$. The result now follows from Lemmas \ref{lemaAB} and \ref{description}.\fim

\begin{remark} Notice that, since $X_A$ is open and compact for each $A\in \mathcal{G}^0$, the characteristic map $1_A$ of the set $X_A$ is an element of $C(X)$. Moreover, each $X_c$, with $c\in \F$, is open and compact (by Proposition \ref{opencompact}), and so the associated characteristic map $1_c$ is also an element of $C(X)$. 
\end{remark}

\begin{lema}\label{densealgebra} The subalgebra $D\subseteq C_0(X)$ generated by all the characteristic maps $1_c$, $1_A$ and $\alpha_c(1_{c^{-1}}1_A)$, with $c\in \bigcup\limits_{n=1}^\infty (\GG^1)^n$ and $A\in\mathcal{G}^0$, is dense in $C_0(X)$. Moreover, for each $0\neq g\in \F$, the subalgebra $D_g$ of $C(X_g)$ generated by all the maps $1_g 1_c$, $1_g 1_A$ and $1_g\alpha_c(1_{c^{-1}}1_A)$, with $c\in \bigcup\limits_{n=1}^\infty (\GG^1)^n$ and $A\in\mathcal{G}^0$, is dense in $C(X_g)$.
\end{lema}

\demo Note that $D$ is a self adjoint subalgebra of $C_0(X)$. By the Stone Weierstrass Theorem it is enough to show that $D$ vanishes nowhere and separates points. 

We show first that $D$ vanishes nowhere. Let $x\in X$. If $|x|=0$ then $x=(A,A)$, where $A\in \mathcal{G}^0$ is minimal infinite emitter. Note that $(A,A)\in X_A$ and that $1_A(x)=1$. If $|x|\geq 1$ then $x=(\beta,B)$. Note that $x\in X_\beta$ and that $1_\beta(x)=1$. So, for each $x\in X$ there exists an element $f\in D$ such that $f(x)\neq 0$, that is, $D$ vanishes nowhere.

Now we show that $D$ separates points. Let $x, y\in X$ with $x\neq y$. Suppose first that $|x|=0$, that is, $x=(A,A)$ where $A\in \mathcal{G}^0$ is a minimal infinite emitter. If $|y|=0$, that is, $y=(B,B)$ where $B\in \mathcal{G}^0$ is a minimal infinite emitter, then $1_A(x)=1$ and $1_A(y)=0$, since $(B,B)\notin X_A$. If $|y|\geq 1$ then $y=(\beta_1...\beta_{|\beta|},B)$ or $y=\beta_1\beta_2...$ where $\beta_i$ are edges. In this case, $1_{\beta_1}(x)=0$ and $1_{\beta_1}(y)=1$. Now suppose that $|x|\geq 1$ and $|y|\geq 1$, with $|x|\leq |y|$. Then $x=(\gamma_1...\gamma_{|\gamma|},A)$ or $x=\gamma_1\gamma_2...$ and $y=(\beta_1...\beta_{|\beta|},A)$ or $y=\beta_1\beta_2...$.
If $\gamma_i\neq \beta_i$ for some $i$ then $1_{\gamma_1...\gamma_i}(x)=1$ and $1_{\gamma_1...\gamma_i}(y)=0$. If $|x|<|y|$ and $\gamma_i=\beta_i$, for $1\leq i \leq |x|$, then $1_{\beta_1...\beta_{|x|+1}}(x)=0$ and $1_{\beta_1...\beta_{|x|+1}}(y)=1$. So we are left with the case $|x|=|y|<\infty$, $x=(\gamma_1...\gamma_{|\gamma|},A)$ and $y=(\gamma_1...\gamma_{|\gamma|},B)$, where $A,B\in \mathcal{G}^0$ are minimal infinite emitters and $\gamma_i$ are edges. Let $\gamma=\gamma_1...\gamma_{|\gamma|}$. Then $$\alpha_\gamma(1_\gamma^{-1} 1_A))(x)=1_{\gamma^{-1}}(\theta_\gamma^{-1}(x))1_A(\theta_\gamma^{-1}(x))=1$$ and 
$$\alpha_\gamma(1_\gamma^{-1} 1_A))(y)=1_{\gamma^{-1}}(\theta_\gamma^{-1}(y))1_A(\theta_\gamma^{-1}(y))=0.$$
So $D$ separates points and it follows that $D$ is dense in $C_0(X)$. 

To see that $D_g$ is dense in $C(X_g)$, with $0\neq g\in \F$, note that for $x\in X_g$ there is a map $f\in D$ such that $f(x)\neq 0$. Hence $1_g(x)f(x)=f(x)\neq 0$. Furthermore, for $x,y\in X_g$, with $x\neq y$, there is a map $h\in D$ such that $h(x)=1$ and $h(y)=0$. Hence $1_g(x)h(x)=1$ and $1_g(y)h(y)=0$. Therefore $D_g$ vanishes nowhere and separates points and hence $D_g$ is dense in $C(X_g)$.\fim

We can now prove the main theorem of this section, which generalizes \cite[Theorem~4.11]{GRUltra}. A version of this theorem for directed graphs (possibly with sinks) can be seen in \cite[Theorem~3.1]{toke}.

\begin{teorema}\label{crossedproduct} Let $\GG$ be an ultragraph with no sinks that satisfies Condition (RFUM). Then there exists a *-isomorphism $\Phi:C^*(\mathcal{G})\rightarrow C_0(X)\rtimes_\alpha \F$ such that $\Phi(s_e)=1_e\delta_e$, for each edge $e\in \mathcal{G}^1$, and $\Phi(p_A)=1_A\delta_0$, for all $A\in \GG^0$.
\end{teorema}

\demo Define $\Phi(s_e)=1_e\delta_e$, for each edge $e\in \mathcal{G}^1$, and $\Phi(p_A)=1_A\delta_0$, for each $A\in \mathcal{G}^0$. By the universality of $C^*(\GG)$, to check that $\Phi$ extends to $C^*(\mathcal{G})$ it is enough to verify that the set $\{\Phi(s_e)\}_{e\in \mathcal{G}^1}\cup \{\Phi(p_A)\}_{A\in \mathcal{G}^0}$ satisfies the conditions of Definition \ref{def of C^*(mathcal{G})}. First note that $\Phi(s_e)\Phi(s_e)^*=1_e\delta_e1_{e^{-1}}\delta_{e^{-1}}=1_e\delta_0$ and $\Phi(s_e)\Phi(s_e)^*\Phi(s_e)= 1_e \delta_e$. Hence $\{\Phi(s_e): e\in \GG^1\}$ is a family of partial isometries with orthogonal ranges. 

The first condition of Definition \ref{def of C^*(mathcal{G})} follows from Lemma \ref{lemaAB}. 

To verify the second condition, let $e\in \mathcal{G}^1$ and note that $$\Phi(s_e)^*\Phi(s_e)=(1_e\delta_e)^*1_e\delta_e=$$ $$=\alpha_{e^{-1}}(\alpha_e(1_{e^{-1}})1_e)\delta_0=1_{e^{-1}}\delta_0=1_{r(e)}\delta_0=\Phi(p_{r(e)}).$$

Moreover, since $X_e\subseteq X_{s(e)}$ then $1_e\delta_0 1_{s(e)}\delta_0=1_e\delta_0$, and so we have that $\Phi(s_e)\Phi(s_e)^*\leq \Phi(s(e))$ and the third condition of Definition \ref{def of C^*(mathcal{G})} follows.

To verify the last condition of Definition \ref{def of C^*(mathcal{G})}, let $v\in G^0$ be such that $0<|s^{-1}(v)|<\infty$. Then $X_v=\bigcup\limits_{e\in s^{-1}(v)}X_e$, and since the projections $\{1_e: e \in \GG^1 \}$ are pairwise orthogonal, then $\sum\limits_{e\in s^{-1}(v)}1_e=1_v$. Hence

$$\sum\limits_{e\in s^{-1}(v)}\Phi(s_e)\Phi(s_e)^*=\sum\limits_{e\in s^{-1}(v)}1_e\delta_0=1_v\delta_0=\Phi(p_v).$$

Next we show that $\Phi$ is surjective. It is not hard to see that for $a=a_1...a_n,b=b_1...b_m\in P$, where $a_i,b_j$ are edges, we get $$\Phi(s_a)=\Phi(s_{a_1})...\Phi(s_{a_n})=1_a\delta_a,$$ 
$$\Phi(s_a)\Phi(s_b)^*=\Phi(s_{a_1})...\Phi(s_{a_n})\Phi(b_m)^*...\Phi(s_{b_1})^*=1_{ab^{-1}}\delta_{ab^{-1}}$$ and $$\Phi(s_a)\Phi(s_b)^*\Phi(s_b)\Phi(s_a)^*=1_{ab^{-1}}\delta_0.$$ So, for all $c\in \F$, $1_c\delta_0$ and $1_c\delta_c$ belong to $\text{Im}(\Phi)$. Hence, for each $c\in \F$ and $A\in \mathcal{G}^0$, we have that $\alpha_c(1_{c^{-1}}1_A)\delta_c=1_c\delta_c1_A\delta_0\in \text{Im}(\Phi)$ and $\alpha_c(1_{c^{-1}}1_A)\delta_0=\alpha_c(1_{c^{-1}}1_A)\delta_c 1_{c^{-1}}\delta_{c^{-1}}\in \text{Im}(\Phi)$. By Lemma $\ref{densealgebra}$ we get that $C_0(X)\delta_0\subseteq \text{Im}(\Phi)$. 

Now note that for $g\in \F$, $g\neq 0$, we have that $(1_g1_c)\delta_c=1_c\delta_0 1_g\delta_g\in \text{Im}(\Phi)$, $(1_g1_A)\delta_g=1_A\delta_0 1_g\delta_g\in \text{Im}(\Phi)$ and $1_g\alpha_c(1_{c^{-1}}1_A)\delta_g=\alpha_c(1_{c^{-1}}1_A)\delta_01_g\delta_g\in \text{Im}(\Phi)$, for all $c\in \F$ and for all $A\in \mathcal{G}^0$. By Lemma \ref{densealgebra} it follows that $C(X_g)\delta_g\subseteq \text{Im}(\Phi)$. So $\Phi$ is surjective. 

Finally we show that $\Phi$ is injective, using the Gauge-Invariant Uniqueness Theorem of \cite{Tom3}. Let $S^1$ be the unit circle and $\varphi:S^1\rightarrow Aut(C^*(\mathcal{G}))$ be the Gauge action. Recall that, for each $z\in S^1$, we have $\varphi_z(s_e)=zs_e$, for each edge $e\in \GG^1$, and $\varphi_z(p_A)=p_A$, for each $A\in\mathcal{G}^0$. To use the Gauge-Invariant Uniqueness Theorem we need to define, for each $z\in S^1$, an automorphism $\psi_z:C_0(X)\rtimes_\alpha \F\rightarrow C_0(X)\rtimes_\alpha \F$. We do this by first defining $\psi_z$ in the dense subalgebra $C_0(X)\rtimes_{alg}\F$ (which consists of finite sums of the form $\sum a_g \delta_g$) and then extending it to the partial crossed product.

For $g=ab^{-1}\in \F$ with $a, b\in P\cup \{0\}$, write $a=a_1...a_n$, $b=b_1...b_m$, where $a_i, b_j$ are edges, and define $|g|=|a_1...a_nb_m^{-1}...b_1^{-1}|=n-m$. For $z\in S^1$ define $\psi_z:C_0(X)\rtimes_{alg}\F\rightarrow C_0(X)\rtimes_\alpha \F$ by $\psi_z(\sum\limits_g a_g\delta_g)=\sum\limits_g z^{|g|}a_g\delta_g$, which is a *-homomorphism. Note that $||\psi_z(\sum\limits_g a_g\delta_g)||=||\sum\limits_g z^{|g|}a_g\delta_g||\leq \sum\limits_g||a_g\delta_g||$. By \cite[Corollary 17.11]{exelbook}, $||a_g\delta_g||=||a_g||$ for each $g\in \F$, and then 
$||\psi_z(\sum\limits_g a_g\delta_g)||\leq\sum\limits_g ||a_g||$. So, $\psi_z$ extends to a *-homomorphism $\psi_z:C_0(X)\rtimes_\alpha \F\rightarrow C_0(X)\rtimes_\alpha \F$. 

To end the proof notice that the map $S^1\ni z\mapsto \psi_z \in Aut(C_0(X)\rtimes_\alpha \F)$ is a strongly continuous action of $S^1$. Moreover, $\Phi\circ \varphi_z=\psi_z\circ \Phi$ for each $z\in S^1$. By \cite[Corollary 17.11]{exelbook} we have that $1_A\delta_0\neq 0$ for each $A\in \mathcal{G}^0$ and hence, by the Gauge-Invariant Uniqueness Theorem of \cite{Tom3}, $\Phi$ is injective. \fim

\begin{remark} To prove injectivity in the theorem above we could have also used the general Cuntz-Krieger Theorem for ultragraphs showed in \cite{dhd}.
\end{remark}

\section{An invariant for shift conjugacy}

In this section we prove that if two ultragraph edge shift spaces are conjugate, via a conjugacy that preserves length, then their associated ultragraph C*-algebras are isomorphic. As a consequence we deduce that there exist ultragraph edge shift spaces that are not conjugate, via a conjugacy that preserves length, to the edge shift space of a graph. Before we prove our main result we need the following lemma.

\begin{lema}\label{conjshiftcomute}  Let $\GG_1, \GG_2$ be two ultragraphs and $X$ and $Y$ be the associated edge shift spaces, respectively. If $\phi: X \rightarrow Y$ is a conjugacy that preserves length then, for all $a \in \GG^1_1$, and for all $x$ such that $|x|\geq 1$ or $|x|=0$ and $x\subseteq r(a)$, there exists $b\in \GG_2^1$ such that $\phi(a \cdot x) = b \cdot \phi (x).$
\end{lema}

\demo Notice that, if $|x|\geq 1$ then $\phi(ax) =  b_1 b_2\ldots$. Since $\phi$ commutes with the shift we have that $\phi(x)=\phi(\sigma(ax)) = \sigma (\phi(ax)) = b_2 \ldots$ and the result follows.

If $|x|=0$ and $x\subseteq r(a)$ then $\phi (a \cdot x) = (b, B)$. Also $\phi(x)= \phi(r(a)\cap x)=\phi(\sigma (a\cdot x)) =  \sigma (\phi(a\cdot x)) = (B,B)$. Now notice that $\phi(a \cdot x)= (b, r(b)) \cdot (B,B) = b \cdot \phi (x)$. \fim

\begin{teorema} Let $\GG_1, \GG_2$ be two ultragraphs with no sinks that satisfy Condition (RFUM), and such that their shift spaces, $X$ and $Y$ respectively, are conjugate via a conjugacy $\phi:X\rightarrow Y$ that preserves length. Then $C^*(\GG_1)$ and $C^*(\GG_2)$ are isomorphic, via an *-isomorphism that intertwines the canonical Gauge actions and maps the commutative C*-subalgebra of $C^*(\GG_1)$, generated by $\{s_{e_1}...s_{e_n}p_As_{e_n}^*...s_{e_1}^*: e_i\in \mathcal{G}_1^1, A\in \mathcal{G}^0\}$, to the corresponding C*-subalgebra of $C^*(\GG_2)$.
\end{teorema}

\demo Let $\F_1$ and $\F_2$ be the free groups generated by $\GG_1^1$ and $\GG_2^1$, respectively. By Theorem \ref{crossedproduct}, $C^*(\GG_1)$ and $C^*(\GG_2)$ are *-isomorphic to the partial crossed products $C_0(X)\rtimes_\beta \F_1$ and $C_0(Y)\rtimes_\alpha \F_2$, respectively.

We will show that the $C^*$-algebras $C^*(\GG_1)$ and $C_0(Y)\rtimes_\alpha \F_2$ are isomorphic. 

First, using the universality of $C^*(\GG_1)$, we show that there exists a *-homomorphism $\pi:C^*(\GG_1)\rightarrow C_0(Y)\rtimes_\alpha \F_2$. 

Note that for each $A\in \GG_1^0$ the set $\phi(X_A)$ is clopen, so that $1_{\phi(X_A)}\in C_0(Y)$. Define $$\pi(p_A)=1_{\phi(X_A)}\delta_0.$$
Fix an $e\in \GG_1^1$. Note that, since $\phi$ preserves length,  $\phi(X_e)\subseteq \bigcup\limits_{a\in \GG_2^1}Y_a$. Since $\phi(X_e)$ is compact, and each $Y_a$ is open, the set $F_e=\{a\in \GG_2^1:\phi(X_e)\cap Y_a\neq \emptyset\}$ is finite. Moreover, the sets $Y_a$ are pairwise disjoint. Define $$\pi(s_e)=\sum_{a\in F_e}1_{\phi(X_e)}1_a\delta_a.$$

To show that $\pi$ extends to a homomorphism we need to show that $\{\pi(p_A):A\in \GG_1^0\}$ and $\{\pi(S_e): e \in \GG^1_1\}$ satisfy the conditions of Definition \ref{def of C^*(mathcal{G})} for $\GG_1$.

For each $A,B\in \GG_1^0$ it holds that $\pi(p_A)\pi(p_B)=1_{\phi(X_A)\cap \phi(X_B)}\delta_0=1_{\phi(X_A\cap X_B)}\delta_0$. By Lemma~$\ref{lemaAB}$, $X_A\cap X_B=X_{A\cap B}$, and then $1_{\phi(X_A\cap X_B)}\delta_0=1_{\phi(X_{A\cap B})}\delta_0=\pi(p_{A\cap B})$.
Moreover, $\pi(p_{A\cup B})=1_{\phi(X_{A\cup B})}\delta_0$ and, since $X_{A\cup B}=X_A\cup X_B$, by Lemma~\ref{lemaAB}, we get $1_{\phi(X_{A\cup B})}\delta_0=1_{\phi(X_A\cup X_B)}\delta_0=1_{\phi(X_A)\cup \phi(X_B)}\delta_0=(1_{\phi(X_A)}+1_{\phi(X_B)}-1_{\phi(X_A)\cap \phi(X_B)})\delta_0=(1_{\phi(X_A)}+1_{\phi(X_B)}-1_{\phi(X_{A\cap B})})\delta_0=\pi(p_A)+\pi(p_B)-\pi(p_{A\cap B})$. So condition $1$ of Definition \ref{def of C^*(mathcal{G})} holds. 

In what comes next we use the following notation: if $Z$ is a set, $R\subseteq Z$, and $x\in Z$, then $[x\in R]=1$ if $x\in R$, and $[x\in R]=0$ else.  

Next, note that $\pi(s_e)^*\pi(s_e)=\sum\limits_{a\in F_e}\alpha_{a^{-1}}(1_a1_{\phi(X_e)})\delta_0$, and $$\sum\limits_{a\in F_e}\alpha_{a^{-1}}(1_a1_{\phi(X_e)})(x)=\sum\limits_{a\in F_e}[x\in Y_{a^{-1}}][\theta_a(x)\in\phi(X_e)]=$$ $$=\sum\limits_{a\in F_e}[x\in Y_{a^{-1}}][ax\in\phi(X_e)]=\sum\limits_{a\in F_e}[x\in Y_{a^{-1}}][\phi^{-1}(ax)\in X_e].$$ Notice that if a term on the right hand side of the equation above is non-zero then $x\in Y_{a^{-1}}$, so if $|x|=0$ then $x\subseteq r(a)$. Therefore we can use Lemma~\ref{conjshiftcomute} to obtain that $\phi^{-1}(ax)=z_a\phi^{-1}(x)$ (since $\phi^{-1}$ is a conjugacy) and so $\phi^{-1}(ax)\in X_e$ if, and only if, $z_a=e$. Hence $\sum\limits_{a\in F_e}[x\in Y_{a^{-1}}][\phi^{-1}(ax)\in X_e]=0$ or 
$\sum\limits_{a\in F_e}[x\in Y_{a^{-1}}][\phi^{-1}(ax)\in X_e]=1$.
We need to show that $\sum\limits_{a\in F_e}[x\in Y_{a^{-1}}][\phi^{-1}(ax)\in X_e]=[x\in \phi(X_{r(e)})]$, for each $x\in Y$.

First suppose that $\sum\limits_{a\in F_e}[x\in Y_{a^{-1}}][\phi^{-1}(ax)\in X_e]=1$. Then $[\phi^{-1}(ax)\in X_e]=1$ for some $a\in F_e$. Hence $z_a\phi^{-1}(x)\in X_e$, that is, $z_a=e$ and so $z_a\phi^{-1}(x)=e\phi^{-1}(x)\in X_e$ Therefore $\phi^{-1}(x)\in X_{r(e)}$, and so $x\in \phi(X_{r(e)})$.

Now suppose that $[x\in \phi(X_{r(e)})]=1$, that is, $x\in \phi(X_{r(e)})$. Then $\phi^{-1}(x)\in X_{r(e)}=X_{e^{-1}}$, and so $e\phi^{-1}(x)\in X_e$ and $\phi(e\phi^{-1}(x))\in \phi(X_e)$. Let $a\in F_e$ be the unique element such that $\phi(e\phi^{-1}(x))\in Y_a$. 
Since $\phi$ is a conjugacy, by Lemma~\ref{conjshiftcomute} (notice that if $|x|=0$ then $\phi^{-1}(x)\subseteq r(e)$), we have that $\phi(e\phi^{-1}(x))=a x$. So, $x\in Y_{a^{-1}}$ and $\phi(ax)\in X_e$, from where we get that $\sum\limits_{a\in F_e}[x\in Y_{a^{-1}}][\phi^{-1}(ax)\in X_e]=1$. 

So we proved that $\pi(s_e)^*\pi(s_e)=\pi(p_{r(e)})$, which is the second condition of Definition \ref{def of C^*(mathcal{G})}.

To verify conditions 3 and 4 of Definition~\ref{def of C^*(mathcal{G})} we need the following claim, the proof of which we omit since it is essentially the same as the proof of Claim 2 in \cite[Theorem~4.12]{GRUltra}.

{\it Claim 1: Let $e\in \GG_1^1$. Then $\pi(s_e)\pi(s_e)^*=1_{\phi(X_e)}\delta_0$, and $\pi(s_e)1_{\phi(X_c)}\pi(s_e)^*=1_{\phi(X_{ec})}\delta_0$ for each $c\in \bigcup\limits_{n=1}^\infty(\GG_1^1)^n$.}

Since $X_e\subseteq X_{s(e)}$ we get that $\pi(s_e)\pi(s_e)^*=1_{\phi(X_e)}\delta_0\leq 1_{\phi(X_{s(e)})}\delta_0=\pi(p_{s(e)})$, and so condition 3 of Definition \ref{def of C^*(mathcal{G})} holds. Moreover, $\{\pi(s_e)\}_{e\in \GG_1^1}$ is a family of partial isometries with orthogonal ranges. To verify condition 4, let $v$ be a vertex of $\GG_1$ such that $0<|s^{-1}(v)|<\infty$. Then $X_v=\bigcup\limits_{e\in s^{-1}(v)}X_e$ and so $\pi(p_v)=1_{\phi(X_v)}\delta_0=\sum\limits_{e\in s^{-1}(v)}1_{\phi(X_e)}\delta_0=\sum\limits_{e\in s^{-1}(v)}\pi(s_e)\pi(s_e)^*$. 

So, we get a *-homomorphism $\pi:C^*(\GG_1)\rightarrow C_0(Y)\rtimes_\alpha\F_2$.

To verify that $\pi$ is surjective, we need to prove first the following result:

{\it Claim 2: For each $c\in \bigcup\limits_{n=1}^\infty(\GG_1^1)^n$, and $A\in \GG_1^0$, we have that $\pi(s_c)\pi(p_A)\pi(s_c)^*=(\beta_c(1_c^{-1}1_{X_A})\circ \phi^{-1})\delta_0$.}

Let $e\in \GG_1^1$ and $A\in \GG_0^1$. Then $$\pi(s_e)\pi(p_A)\pi(s_e)^*=\sum\limits_{a,b\in F_e}\alpha_a(\alpha_{a^{-1}}(1_{\phi(X_e)}1_a)1_{\phi(X_A)}\alpha_{b^{-1}}(1_{\phi(X_e)}1_b))\delta_{ab^{-1}}.$$

Let $x\in Y_{ab^{-1}}$ and write $x=ay$, with $y\subseteq r(a) \cap r(b)$ if $|y|=0$). Then $$\left(\alpha_a(\alpha_{a^{-1}}(1_{\phi(X_e)}1_a)1_{\phi(X_A)}\alpha_{b^{-1}}(1_{\phi(X_e)}1_b))\right)(x)=$$
$$=(\alpha_{a^{-1}}(1_{\phi(X_e)}1_a)1_{\phi(X_A)}\alpha_{b^{-1}}(1_{\phi(X_e)}1_b))(y)=$$
$$=(1_{\phi(X_e)}(ay)1_{\phi(X_A)}(y)(1_{\phi(X_e)}(by))=$$
$$=[\phi^{-1}(ay)\in X_e][y\in \phi(X_A)][\phi^{-1}(by)\in X_e].$$
Note that, by Lemma~\ref{conjshiftcomute}, $\phi^{-1}(ay)=u\phi^{-1}(y)$ and $\phi^{-1}(by)=v\phi^{
-1}(y)$ (for some $u,v\in \GG_1^1$). So $\phi^{-1}(ay)\in X_e$ if, and only if, $u=e$, and $\phi^{-1}(by)\in X_e$ if, and only if, $v=e$. Hence if $a\neq b$ we get that $$[\phi^{-1}(ay)\in X_e][y\in \phi(X_A)][\phi^{-1}(by)\in X_e]=0,$$
and so $$\pi(s_e)\pi(p_A)\pi(s_e)^*=\sum\limits_{a\in F_e}\alpha_a(\alpha_{a^{-1}}(1_{\phi(X_e)}1_a)1_{\phi(X_A)})\delta_0.$$

Let $x\in \phi(X_e)$. Then $x\in Y_a$ for some $a\in F_e$.
So, $$\pi(s_e)\pi(p_A)\pi(s_e)^*(x)=\alpha_a(\alpha_{a^{-1}}(1_{\phi(X_e)}1_a)1_{\phi(X_A)})(x)=[\sigma(x)\in \phi(X_A)]=$$
$$=[\phi^{-1}(\sigma(x))\in X_A]=[\sigma(\phi^{-1}(x))\in X_A]=\beta_e(1_{e^{-1}}1_{X_A})(\phi^{-1}(x)).$$ For $x\notin \phi(X_e)$ we get $$\pi(s_e)\pi(p_A)\pi(s_e)^*(x)=0=\beta_e(1_{e^{-1}}1_{X_A})(\phi^{-1}(x)).$$ Therefore $\pi(s_e)\pi(p_A)\pi(s_e)^*=(\beta_e(1_{e^{-1}}1_{X_A})\circ\phi^{-1})\delta_0$. Now we proceed by inductive arguments. Suppose that $\pi(s_c)\pi(p_A)\pi(s_c)^*=(\beta_c(1_{c^{-1}}1_{X_A})\circ \phi^{-1})\delta_0$ for each $c\in (\GG_1^1)^n$. Let $b\in \GG_1^1)^{n+1}$ and write $b=ec$, where $|c|=n$ and $e\in \GG_1^1$.
Then $$\pi(b)\pi(p_A)\pi(b)^*=\pi(s_e)\pi(s_c)\pi(p_A)\pi(s_c)^*\pi(s_e)=$$ $$=\pi(s_e)(\beta_c(1_{c^{-1}}1_{X_A})\circ\phi^{-1})\delta_0\pi(s_e)^*=\sum\limits_{d\in F_e}\alpha_d(\alpha_{d^{-1}}(1_{\phi(X_e)}1_d)\beta_c(1_{c^{-1}}1_{X_A})\circ\phi^{-1})\delta_0.$$ Let $x\in \phi(X_e)$ and write $x=ay\in Y_a$, for some $a\in F_e$ (and $y\subseteq r(a)$ if $|y|=0$). Then $$\sum\limits_{d\in F_e}\alpha_d(\alpha_{d^{-1}}(1_{\phi(X_e)}1_d)\beta_c(1_{c^{-1}}1_{X_A})\circ\phi^{-1})(x)=$$ $$=\alpha_a(\alpha_{a^{-1}}(1_{\phi(X_e)}1_a)\beta_c(1_{c^{-1}}1_{X_A})\circ\phi^{-1})(x)=$$
$$=[\sigma^n(\phi^{-1}(\sigma(x)))\in X_{c^{-1}}\cap X_A]=[\sigma^{n+1}(\phi^{-1}(x))\in X_{c^{-1}}\cap X_A]=$$ $$=[\sigma^{n+1}(\phi^{-1}(x))\in X_{b^{-1}}\cap X_A]=\beta_b(1_{b^{-1}}1_{X_A})(\phi^{-1}(x)).$$
For $x\notin \phi(X_e)$ we get that $$\sum\limits_{a\in F_e}\alpha_a(\alpha_{a^{-1}}(1_{\phi(X_e)}1_a)\beta_c(1_{c^{-1}}1_{X_A})\circ\phi^{-1})(x)=0=\beta_b(1_{b^{-1}}1_{X_A})(\phi^{-1}(x)).$$ Therefore $\pi(b)\pi(p_A)\pi(b)^*=\beta_b(1_b^{-1}1_A)\circ \phi^{-1}\delta_0$, and Claim 2 is proved.

We can now show that $\pi$ is surjective. By Lemma~\ref{densealgebra}, the algebra generated by $1_A$, $1_c$ and $\beta_c(1_{c^{-1}}1_A)$, with $A\in \GG_1^0$ and $c\in \bigcup\limits_{n=1}^\infty (\GG_1^1)^n$, is dense in $C_0(X)$. So the algebra generated by $1_{\phi(X_A)}$, $1_{\phi(X_c)}$ and $\beta_c(1_c^{-1}1_{X_A})\circ \phi^{-1}$ is dense in $C_0(Y)$. Notice that, by Claim 1 for $c=c_1...c_n\in \GG_1^n$, we get that $\pi(s_{c_1})...\pi(s_{c_n})\pi(s_{c_n})^*...\pi(s_{c_1})^*=1_{\phi(X_c)}\delta_0$. So $1_{\phi(X_A)}\delta_0$, and $1_{\phi(X_c)}\delta_0$ belong to $\text{Im}(\pi)$, for each $A\in \GG_1^0$ and $c\in \bigcup\limits_{n=1}^\infty (\GG_1^1)^n$. Also, by Claim 2, $\beta_c(1_{c^{-1}}1_{X_A})\delta_0\in \text{Im}(\pi)$. Therefore $C_0(Y)\delta_0\subseteq \text{Im}(\pi)$. Furthermore, notice that $\pi$ maps the C*-subalgebra generated by $\{s_{e_1}...s_{e_n}p_As_{e_n}^*...s_{e_1}^*: e_i\in \mathcal{G}_1^1, A\in \mathcal{G}^0\}$ onto $C_0(Y)\delta_0$.

 In particular, $1_{Y_a}\delta_0\in \text{Im}(\phi)$ for each $a\in \GG_2^1$. Since $Y_a$ is compact, and $\{\phi(X_e)\}_{e\in \GG_1^1}$ is an open cover of $Y_a$, there exists a finite cover $\{\phi(X_e)\}_{e\in H_a}$ of $Y_a$. Then $$\sum\limits_{e\in H_a}1_a\delta_0\pi(s_e)=\sum\limits_{e\in H_a}\sum\limits_{b\in F_e} 1_a1_b1_{\phi(X_e)}\delta_b=\sum\limits_{e\in H_a} 1_a1_{\phi(X_e)}\delta_a=1_a\delta_a,$$ and so $1_a\delta_a\in \text{Im}(\pi)$ for each $a\in \GG_2^1$. Moreover, since $C_0(Y)\delta_0\subseteq \text{Im}(\pi)$ then $1_A\delta_0\in \text{Im}(\pi)$ for each $A\in \GG_2^0$. By Theorem \ref{crossedproduct}, $C_0(Y)\rtimes_\alpha\F_2$ is generated by $1_a\delta_a$, with $a\in \GG_2^1$, and $1_A\delta_0$, with $A\in \GG_2^0$. Hence $\pi$ is surjective.

It remains to show that $\pi$ is injective. Following the arguments at the end of the proof of Theorem \ref{crossedproduct}, we get a group action $\psi:S^1\rightarrow Aut(C_0(Y)\rtimes_\alpha \F_2)$ such that $\psi_z(f_a\delta_a)=zf_a\delta_a$ for each $a\in \GG_1^1$, $f_a\in C(Y_a)$, and $z\in S^1$, and $\psi_z(f_0\delta_0)=f_0\delta_0$ for each $f_0\in C_0(Y)$. Let $\varphi:S^1\rightarrow Aut(C^*(\GG_1))$ be the Gauge action. Note that $\pi$ intertwines the Gauge actions $\varphi$ and $\psi$, and the injectivity of $\pi$ now follows from the Gauge-Invariant Uniqueness Theorem of \cite{Tom3}. \fim

\begin{corolario}\label{newshift} There exist ultragraph edge shift spaces, associated to ultragraphs with no sinks that satisfy Condition (RFUM), that are not conjugate, via a conjugacy that preserves length, to the path space of a graph (viewed as a shift space).
\end{corolario}

\demo Consider the shift space associated to the ultragraph whose C*-algebra is neither an Exel-Laca nor a graph algebra (see \cite{TomSimple}). If this shift space is conjugate, via a conjugacy that preserves length, to the path space of a graph then, by the theorem above, the associated C*-algebras would be isomorphic, which contradicts \cite[Corollary~5.5]{TomSimple}. \fim

Daniel Gon\c{c}alves, Departamento de Matem\'atica, Universidade Federal de Santa Catarina, Florian\'opolis, 88040-900, Brasil

Email: daemig@gmail.com

\vspace{0.5pc}
Danilo Royer,  Departamento de Matem\'atica, Universidade Federal de Santa Catarina, Florian\'opolis, 88040-900, Brasil

Email: danilo.royer@ufsc.br
\vspace{0.5pc}

\end{document}